# The Identification Capacity and Resolvability of Channels with Input Cost Constraint [†]

Te Sun HAN[‡]

May 8, 2000

---

**Abstract:** Given a general channel, we first formulate the idetification capacity problem as well as the resolvability problem with input cost constraint in as the general form as possible, along with relevant fundamental theorems. Next, we establish some mild sufficient condition for the key lemma linking the identification capacity with the resolvability to hold for the continuous input alphabet case with input cost constraint. Under this mild condition, it is shown that we can reach the *continuous*-input fundamental theorem of the same form as that for the fundamental theorem with *finite* input alphabet. Finally, as important examples of this continuous-input fundamental theorem, we show that the identification capacity as well as the resolvability coincides with the channel capacity for stationary additive white (and also non-white) Gaussian noise channels.

# 1 Introduction

In 1989 Ahlswede and Dueck [1] have first devised the novel concept of *identification code* for channels, which is completely different from the traditional transmission code of Shannon [2]. Surprizingly enough, although by slightly relaxing the criterion for reliable transmission, they have shown that a channel with identification codes is capable of reliably transmitting the *doubly* exponential number of messages in blocklength $n$. This is in sharp contrast with the fact that a channel with transmission codes is usually capable of reliably transmitting only the exponential number of messages in blocklength $n$. Intuitively spaeking, the identification code is a system which simultaneously superposes a large number of *simple* hypothesis testings over a channel, whereas the transmission code is a system which conveys the only one *multiple* hypothesis testing over a channel.

On the other hand, Han and Verdú [4] have first coined the concept of *resolvability* for channels. The channel resolvability is defined as the smallest number (per channel input letter) of fair coin flips to generate a channel input of blocklength $n$ whose channel output well approximates any realizable channel output under suitable approximation measures such as the variational distance. A primitive germ of channel resolvability is found in Han and Verdú [3].

Although the resolvability problem looks seemingly very different from the identification problem, we may regard one as the *dual* of the other. In fact, as was shown in Han and Verdú [4], the direct theorem for the channel resolvability gives the converse theorem for the identification capacity, and, on the other hand, the direct theorem for the identification capacity gives the converse theorem for the channel resolvability. However, we should emphasize that this nice duality heavily relies upon the *finiteness* of the channel input alphabet.

In the present paper, we attempt to dispense with this finiteness assumption. In particular, as a typical illustrative case, we consider the general channl with *continuous* input alphabet. Usually, in the case of continuous input alphabet, we need to introduce the input cost constraint such as input power constraint, because, otherwise, the channel capacity as well as the identification capacity and the resolvability may diverge to $\infty$.

For this reason, in Sections 2∼ 5, we first formulate the idetification capacity problem as well as the resolvability problem with input cost constraint in the as general form as possible, along with relevant fundamental theorems. All the definitions along with all the theorems stated in Sections $2 \sim 5$ are straightforward generalizations of those developed already in Han and Verdú [4] for the case *without* input cost constraint. Therefore, our main point is not in this part.

Next, in Section 6, we establish some mild sufficient condition for the key lemma linking the identification capacity with the resolvability to hold for the continuous input alphabet case with input cost constraint. This part is exactly the main purpose of the present paper. Under this mild condition, it is shown that we can reach the *continuous*-input fundamental theorem of the same form as that for the fundamental theorem with *finite* input alphabet. Finally, as important examples of this continuous-input fundamental theorem, we show that the identification capacity as well as the resolvability coincides with the channel capacity for stationary additive white (and also *non*-white) Gaussian noise channels.

## 2 Channel capacity with input cost

In this paper, we consider a very wide class of channels as follows. Let input alphabet $\mathcal{X}$ and output alphabet $\mathcal{Y}$ be arbitrary *abstract* sets. A *general channel* $\mathbf{W} = \{W^n\}_{n=1}^{\infty}$ is a sequence of $n$-dimensional transition probability matrices $W^n : \mathcal{X}^n \to \mathcal{Y}^n$ such that

$$\sum_{\mathbf{y} \in \mathcal{Y}^n} W^n(\mathbf{y}|\mathbf{x}) = 1 \quad (\forall \mathbf{x} \in \mathcal{X}^n).^* \tag{2.1}$$

The class of channels thus defined includes all nonstationary and/or nonergodic channels with arbitrary memory structures (cf. Han and Verdú [4]).

Let $\mathcal{M}_n = \{1, 2, \cdots, M_n\}$ be a message set for this channel, and let $\varphi_n : \mathcal{M}_n \to \mathcal{X}^n$ and $\psi_n : \mathcal{Y}^n \to \mathcal{M}_n$ be the *encoder* and the *decoder*, respectively. For each $i \in \mathcal{M}_n$, the set $\mathcal{D}_i \equiv \psi_n^{-1}(i) \subset \mathcal{Y}^n$ is called the *decoding set* for the message $i$. The probability of error $\varepsilon_n$ with the pair $(\varphi_n, \psi_n)$ is defined by

$$\varepsilon_n = \frac{1}{M_n} \sum_{i=1}^{M_n} W^n(\mathcal{D}_i^c | \varphi_n(i)), \tag{2.2}$$

where "c" denotes the complement of a set. A pair $(\varphi_n, \psi_n)$ with message set $\mathcal{M}_n = \{1, 2, \cdots, M_n\}$ of size $M_n$ and probability of error $\varepsilon_n$ is called an $(n, M_n, \varepsilon_n)$ code for the channel $\mathbf{W} = \{W^n\}_{n=1}^{\infty}$.

---

*In the case where the output alphabet $\mathcal{Y}$ is *abstract*, $W^n(\mathbf{y}|\mathbf{x})$ is understood to be the (conditional) probability measure element $W^n(d\mathbf{y}|\mathbf{x})$ that is measurable in $\mathbf{x}$. Accordingly, the corresponding summation $\sum$ is understood to be the integral $\int$.

Here we introduce the general cost function as follows. Let $c_n : \mathcal{X}^n \to \mathbf{R}$ ($n = 1, 2, \cdots$) be arbitrary functions and set $\mathbf{c} = \{c_n\}_{n=1}^\infty$, which we call the *general* cost function, where $\mathbf{R}$ denotes the set of real numbers. For an $\mathbf{x} \in \mathcal{X}^n$, $c_n(\mathbf{x})$ is called the cost of $\mathbf{x}$, and $\frac{1}{n}c_n(\mathbf{x})$ is called the cost of $\mathbf{x}$ per input letter. The $(n, M_n, \varepsilon_n)$ code that satisfies the *input cost constraint*:

$$\frac{1}{n}c_n\left(\varphi_n(i)\right) \leq \Gamma \quad (\forall n = 1, 2, \cdots; \ \forall i \in \mathcal{M}_n) \tag{2.3}$$

is called an $(n, M_n, \varepsilon_n, \Gamma)$ code, where $\Gamma$ is an arbitrary prescribed constant.

Let us now define the $(\varepsilon, \Gamma)$-channel capacity $C_s(\varepsilon, \Gamma|\mathbf{W})$ of the channel $\mathbf{W} = \{W^n\}_{n=1}^\infty$ *with input cost* as follows, where $0 \leq \varepsilon < 1$ is an arbitrary prescribed constant.

**Definition 2.1**

$R$ is $(\varepsilon, \Gamma)$-achievable $\overset{\text{def}}{\iff}$ There exists an $(n, M_n, \varepsilon_n, \Gamma)$ code
such that $\limsup_{n \to \infty} \varepsilon_n \leq \varepsilon$
and $\liminf_{n \to \infty} \frac{1}{n} \log M_n \geq R.$

**Definition 2.2 ($(\varepsilon, \Gamma)$-channel capacity)**

$$C_s(\varepsilon, \Gamma|\mathbf{W}) = \sup\{R \mid R \text{ is } (\varepsilon, \Gamma)\text{-achievable}\}.$$

In order to demonstrate a general formula for the $(\varepsilon, \Gamma)$-channel capacity, set

$$\mathcal{X}^n(\Gamma) = \left\{\mathbf{x} \in \mathcal{X}^n \ \middle| \ \frac{1}{n}c_n(\mathbf{x}) \leq \Gamma \right\} \tag{2.4}$$

and denote by $\mathcal{S}_\Gamma$ the set of all input processes $\mathbf{X} = \{X^n\}_{n=1}^\infty$ satisfying $\Pr\{X^n \in \mathcal{X}^n(\Gamma)\} = 1$ ($\forall n = 1, 2, \cdots$). Moreover, for any input process $\mathbf{X} = \{X^n\}_{n=1}^\infty$ ($X^n \equiv (X_1^{(n)}, X_2^{(n)}, \cdots, X_n^{(n)})$ takes values in $\mathcal{X}^n$), define

$$J(R|\mathbf{X}) = \limsup_{n \to \infty} \Pr\left\{\frac{1}{n}\log\frac{W^n(Y^n|X^n)}{P_{Y^n}(Y^n)} \leq R\right\}, \tag{2.5}$$

where $Y^n$ is the channel output via $W^n$ due to the channel input $X^n$. [†] Here and hereafter, we use the convention that $P_Z(\cdot)$ denotes the probability distribution of a random variable $Z$. Then, we have the following theorem:

---

[†]In the case where the input and output alphabets $\mathcal{X}, \mathcal{Y}$ are *abstract* (*not* necessarily countable), $\frac{W^n(Y^n|X^n)}{P_{Y^n}(Y^n)}$ in (2.5) is understood to be $g(Y^n|X^n)$, where $g(\mathbf{y}|\mathbf{x}) \equiv \frac{W^n(d\mathbf{y}|\mathbf{x})}{P_{Y^n}(d\mathbf{y})}$ $= \frac{W^n(d\mathbf{y}|\mathbf{x})P_{X^n}(d\mathbf{x})}{P_{Y^n}(d\mathbf{y})P_{X^n}(d\mathbf{x})} = \frac{P_{X^n Y^n}(d\mathbf{x}, d\mathbf{y})}{P_{X^n}(d\mathbf{x})P_{Y^n}(d\mathbf{y})}$ is the Radon-Nikodym derivative that is measurable in $(\mathbf{x}, \mathbf{y})$.

**Theorem 2.1** (Verdú and Han [5], Han [6])  The $(\varepsilon, \Gamma)$-channel capacity $C_s(\varepsilon, \Gamma | \mathbf{W})$ of the channel $\mathbf{W}$ is given by

$$C_s(\varepsilon, \Gamma | \mathbf{W}) = \sup_{\mathbf{X} \in \mathcal{S}_\Gamma} \sup \{R \mid J(R|\mathbf{X}) \leq \varepsilon\} \quad (0 \leq \forall \varepsilon < 1) \qquad (2.6)$$

*Proof:*

The proof basically parallels that of Theorem 6 in Verdú and Han [5], except for that, here, we have to take account of the cost constraint. □

**Remark 2.1** Since a special case where $c_n(\mathbf{x}) = n$ $(\forall \mathbf{x} \in \mathcal{X}^n)$ and $\Gamma = 1$ implies actually *no* input cost constaraint, in this case $C_s(\varepsilon, \Gamma | \mathbf{W})$ reduces to the usual $\varepsilon$-capacity $C(\varepsilon | \mathbf{W})$ *without* input cost constraint. On the other hand, if we consider another special case with $\varepsilon = 0$ and put $C_s(\Gamma | \mathbf{W}) \equiv C_s(0, \Gamma | \mathbf{W})$ $(\varepsilon = 0)$, then $C_s(\Gamma | \mathbf{W})$ is called the $\Gamma$-cost capacity of the channel $\mathbf{W}$, which is defined as the capacity with input cost constraint (2.3) and asymptotically vanishing probability of error ($\lim_{n \to \infty} \varepsilon_n = 0$). It is obvious that

$$C_s(\Gamma | \mathbf{W}) \leq C_s(\varepsilon, \Gamma | \mathbf{W}) \quad (0 \leq \forall \varepsilon < 1). \qquad (2.7)$$

Let us now give a general formula for the $\Gamma$-cost capacity $C_s(\Gamma | \mathbf{W})$. Let $\mathbf{X} = \{X^n\}_{n=1}^{\infty}$ be an arbitrary input process, and define[‡]

$$\underline{I}(\mathbf{X}; \mathbf{Y}) \equiv \text{p-}\liminf_{n \to \infty} \frac{1}{n} \log \frac{W^n(Y^n | X^n)}{P_{Y^n}(Y^n)} \qquad (2.8)$$

where $Y^n$ is the channel output via $W^n$ due to the channel input $X^n$ and we have put $\mathbf{Y} = \{Y^n\}_{n=1}^{\infty}$. Then, as a collorary of Theorem 2.1, we have

**Theorem 2.2** The $\Gamma$-cost capacity $C_s(\Gamma | \mathbf{W})$ of the channel $\mathbf{W}$ is given by

$$C_s(\Gamma | \mathbf{W}) = \sup_{\mathbf{X} \in \mathcal{S}_\Gamma} \underline{I}(\mathbf{X}; \mathbf{Y}). \qquad (2.9)$$

*Proof:* It suffices only to set $\varepsilon = 0$ in (2.6) of Theorem 2.1. □

A general channel $\mathbf{W} = \{W^n\}_{n=1}^{\infty}$ is said to satisfy the *strong converse property* with input cost constraint $\Gamma$ if the probability $\varepsilon_n$ of error for channel coding with any rate $R$ such that $R > C_s(\Gamma|\mathbf{W})$ necessarily approaches one as $n$ tends to $\infty$ (cf. Verdú and Han [5]). In this connection, the following theorem will be used later.

---

[‡]For any sequence $\{Z_n\}_{n=1}^{\infty}$ of real-valued random variables, we define the *limit inferior in probability* (cf. Han and Verdú [4]) of $\{Z_n\}_{n=1}^{\infty}$ by $\text{p-}\liminf_{n \to \infty} Z_n = \sup\{\alpha | \lim_{n \to \infty} \Pr\{Z_n < \alpha\} = 0\}$.

**Theorem 2.3** The necessary and sufficient condition for the channel $\mathbf{W}$ to satisfy the strong converse property with input cost constraint $\Gamma$ is

$$\sup_{\mathbf{X}\in\mathcal{S}_\Gamma} \underline{I}(\mathbf{X};\mathbf{Y}) = \sup_{\mathbf{X}\in\mathcal{S}_\Gamma} \overline{I}(\mathbf{X};\mathbf{Y}), \tag{2.10}$$

where[§]

$$\overline{I}(\mathbf{X};\mathbf{Y}) \equiv \text{p-}\limsup_{n\to\infty} \frac{1}{n}\log\frac{W^n(Y^n|X^n)}{P_{Y^n}(Y^n)}. \tag{2.11}$$

*Proof:*

The proof basically parallels that of Theorem 7 of Verdú and Han [5], except for that, here, we have to take account the cost constraint. □

## 3 Identification capacity with input cost

Given a general channel $\mathbf{W} = \{W^n\}_{n=1}^\infty$ with arbitrary *abstract* input and output alphabets $\mathcal{X}$, $\mathcal{Y}$, let us define the identification code for the channel $\mathbf{W}$ as follows. Let $\mathcal{N}_n = \{1,2,\cdots,N_n\}$ be a message set and let $\mathcal{P}_\Gamma(\mathcal{X}^n)$ denote the set of all probability distributions $Q$ on $\mathcal{X}^n$ that satisfy the cost constraint (cf. (2.4)):

$$Q(\mathcal{X}^n(\Gamma)) = 1. \tag{3.1}$$

We first choose $N_n$ probability distributions $Q_1, Q_2, \cdots, Q_{N_n} \in \mathcal{P}_\Gamma(\mathcal{X}^n)$; and when we want to send message $i \in \mathcal{N}_n$, the (*stochastic*) encoder $\varphi_n : \mathcal{N}_n \to \mathcal{X}^n$ generate at random an input sequence $\mathbf{x} \in \mathcal{X}^n$ according to the probability distribution $Q_i$: $Q_i = \varphi_n(i)$, where $Q_i = \varphi_n(i)$ is called the *codeword* for the message $i$ and $\mathcal{C}_n = \{Q_1, Q_2, \cdots, Q_{N_n}\}$ is called the *code*. On the other hand, at the side of the receiver, we have $N_n$ (deterministic) decoders $\psi_n^{(i)}$ ($i = 1, 2, \cdots, N_n$) and $N_n$ decoding sets $\mathcal{D}_i \subset \mathcal{Y}^n$ ($i = 1, 2, \cdots, N_n$), where $\mathcal{D}_1, \mathcal{D}_2, \cdots, \mathcal{D}_{N_n}$ are *not* required to be mutually disjoint. Each decoder $\psi_n^{(i)}$ ($i = 1, 2, \cdots, N_n$) works as follows: Let $\mathbf{y} \in \mathcal{Y}^n$ be the received output sequence. If $\mathbf{y} \in \mathcal{D}_i$, then the decoder $\psi_n^{(i)}$ declares that the message $i$ was sent. Otherwise, the decoder $\psi_n^{(i)}$ declares that the message $i$ was *not* sent. It should be noted here that each decoder $\psi_n^{(i)}$ is interested in whether the corresponding message $i$ was sent or not, but not interested in what message was sent if the decoder $\psi_n^{(i)}$ declares that the message $i$ was not sent.

---

[§]For any sequence $\{Z_n\}_{n=1}^\infty$ of real-valued random variables, we define the *limit superior in probability* (cf. Han and Verdú [4]) of $\{Z_n\}_{n=1}^\infty$ by $\text{p-}\limsup_{n\to\infty} Z_n = \inf\{\beta | \lim_{n\to\infty}\Pr\{Z_n > \beta\} = 0\}$.

Here, for notational simplicity, we use the convention that, for a probability distribution $Q \in \mathcal{P}_\Gamma(\mathcal{X}^n)$, $QW^n$ denotes the probability distribution on $\mathcal{Y}^n$ such that

$$QW^n(\mathbf{y}) \equiv \sum_{\mathbf{x} \in \mathcal{X}^n} Q(\mathbf{x}) W^n(\mathbf{y}|\mathbf{x}) \quad (\forall \mathbf{y} \in \mathcal{Y}^n). \tag{3.2}$$

We now define

$$\mu_n^{(i)} = Q_i W^n(\mathcal{D}_i^c) \quad (i = 1, 2, \cdots, N_n), \tag{3.3}$$
$$\lambda_n^{(j,i)} = Q_j W^n(\mathcal{D}_i) \quad (j \neq i), \tag{3.4}$$

and set

$$\mu_n = \max_{1 \leq i \leq N_n} \mu_n^{(i)}, \tag{3.5}$$
$$\lambda_n = \max_{1 \leq i \neq j \leq N_n} \lambda_n^{(i,j)}. \tag{3.6}$$

We call $\mu_n, \lambda_n$ the first kind of error probability and the second kind of error probability, respectively, and call the pair $(\varphi_n, \psi_n)$ with error probabilities $\mu_n, \lambda_n$ an $(n, N_n, \mu_n, \lambda_n, \Gamma)$ *identification code*, where we have put

$$\psi_n \equiv (\psi_n^{(1)}, \psi_n^{(2)}, \cdots, \psi_n^{(N_n)}).$$

Usually, we impose on the error brobabilities $\mu_n, \lambda_n$ the condition of the following form:

$$\limsup_{n \to \infty} \mu_n \leq \mu, \tag{3.7}$$
$$\limsup_{n \to \infty} \lambda_n \leq \lambda, \tag{3.8}$$

where $0 \leq \mu < 1$, $0 \leq \lambda < 1$ are any prescribed constants. We now define the $(\mu, \lambda, \Gamma)$-identification capacity $D(\mu, \lambda, \Gamma|\mathbf{W})$ of the channel $\mathbf{W}$ as follows.

**Definition 3.1**

$R$ is $(\mu, \lambda, \Gamma)$-achievable $\overset{\text{def}}{\iff}$ There exists an $(n, N_n, \mu_n, \lambda_n, \Gamma)$ code
such that $\limsup_{n \to \infty} \mu_n \leq \mu$, $\limsup_{n \to \infty} \lambda_n \leq \lambda$
and $\liminf_{n \to \infty} \frac{1}{n} \log \log N_n \geq R.$

**Definition 3.2 ($(\mu, \lambda, \Gamma)$-identification capacity)**

$$D(\mu, \lambda, \Gamma | \mathbf{W}) = \sup \{R \mid R \text{ is } (\mu, \lambda, \Gamma)\text{-achievable}\}.$$

With these preparations, we have the following theorem connecting the identification capacity with the channel capacity:

**Theorem 3.1 (Direct theorem)** Let $\mathbf{W}$ be any general channel. If $0 \leq \varepsilon \leq \mu$ and $0 \leq \varepsilon \leq \lambda$, then it holds that

$$D(\mu, \lambda, \Gamma | \mathbf{W}) \geq C_s(\varepsilon, \Gamma | \mathbf{W}). \qquad (3.9)$$

*Proof:*

The proof basically parallels that of Theorem 1 in Ahlswede and Dueck [1], except for that, here, we have to take account of the cost constraint (also, cf. Han and Verdú [4]).  □

## 4 Channel resolvability with input cost

In this section, let us formulate the *resolvability* problem via channels with input cost constraint $\Gamma$. Suppose that we are given a general channel $\mathbf{W} = \{W^n\}_{n=1}^\infty$ with arbitrary *abstract* input and output alphabets $\mathcal{X}$, $\mathcal{Y}$, and let $\mathbf{X} = \{X^n\}_{n=1}^\infty \in \mathcal{S}_\Gamma$ be an arbitrarily given input process, that is, $\Pr\{X^n \in \mathcal{X}^n(\Gamma)\} = 1$ for all $n = 1, 2, \cdots$ (input cost constraint $\Gamma$), where $\mathcal{X}^n(\Gamma)$ is as specified in (2.4). Denote by $\mathbf{Y} = \{Y^n\}_{n=1}^\infty$ the output process via $\mathbf{W} = \{W^n\}_{n=1}^\infty$ due to the input process $\mathbf{X} = \{X^n\}_{n=1}^\infty$. We want to approximate the distribution of the output $Y^n$ via $W^n$ due to the channel input $X^n$ by the distribution of the output via $W^n$ due to another appropriate channel input as follows: Let $U_{M_n}$ be the random variable that is *uniformly distributed* on the set $\mathcal{M}_n \equiv \{1, 2, \cdots, M_n\}$. We call this $U_{M_n}$ the *uniform random number* of size $M_n$. The encoder $\varphi_n : \mathcal{M}_n \to \mathcal{X}^n(\Gamma)$ deterministically maps $U_{M_n}$ into the channel input $\tilde{X}^n \equiv \varphi_n(U_{M_n})$, where $\Pr\{\tilde{X}^n \in \mathcal{X}^n(\Gamma)\} = 1$. Denote by $\tilde{Y}^n$ the channel output via $W^n$ due to the channel input $\tilde{X}^n$, and set $\tilde{\mathbf{X}} = \{\tilde{X}^n\}_{n=1}^\infty$, $\tilde{\mathbf{Y}} = \{\tilde{Y}^n\}_{n=1}^\infty$. It is evident that $\tilde{\mathbf{X}} = \{\tilde{X}^n\}_{n=1}^\infty \in \mathcal{S}_\Gamma$. Usually, we require

$$\lim_{n\to\infty} d(Y^n, \tilde{Y}^n) = 0, \tag{4.1}$$

where $d(Z, V)$ denotes the *variational distance* between the probability distributions $P_Z, P_V$ of random variables $Z, V$. More generally, condition (4.1) can be relaxed to

$$\limsup_{n\to\infty} d(Y^n, \tilde{Y}^n) \leq \delta, \tag{4.2}$$

where $0 \leq \delta < 2$ is an arbitrary prescribed constant. Under condition (4.2) we want to make the size $M_n$ of the uniform random number $U_{M_n}$ as small as possible. Thus, we formulate the resolvability problem as follows.

**Definition 4.1**

$R$ is $(\delta, \Gamma)$-achievable for input $\mathbf{X} = \{X^n\}_{n=1}^\infty \in \mathcal{S}_\Gamma$ $\overset{\text{def}}{\iff}$ There exists an encoder $\varphi_n$ such that $\Pr\{\tilde{X}^n \equiv \varphi_n(U_{M_n}) \in \mathcal{X}^n(\Gamma)\} = 1$,
$$\limsup_{n\to\infty} d(Y^n, \tilde{Y}^n) \leq \delta \text{ and}$$
$$\limsup_{n\to\infty} \frac{1}{n} \log M_n \leq R,$$

where $Y^n$, $\tilde{Y}^n$ are the channel outputs due to the channel inputs $X^n$, $\tilde{X}^n$, respectively.

**Definition 4.2**

$R$ is $(\delta, \Gamma)$-achievable $\overset{\text{def}}{\iff}$ $R$ is $(\delta, \Gamma)$-achievable for all inputs $\mathbf{X} \in \mathcal{S}_\Gamma$.

**Definition 4.3 (($\boldsymbol{\delta, \Gamma}$)-channel resolvability)**

$$S(\delta, \Gamma | \mathbf{W}) = \inf\{R \mid R \text{ is } (\delta, \Gamma)\text{-achievable}\}.$$

With these definitions, we have the following theorem.

**Theorem 4.1 (Direct theorem)** Let $\mathbf{W}$ be any general channel. Then, for all $\delta \geq 0$ we have

$$S(\delta, \Gamma | \mathbf{W}) \leq \sup_{\mathbf{X} \in \mathcal{S}_\Gamma} \overline{I}(\mathbf{X}; \mathbf{Y}). \tag{4.3}$$

*Proof:*

The proof basically parallels that of Theorem 4 in Han and Verdú [4], except for that, here, we have to take account of the cost constraint. □

# 5 Fundamental theorem

So far, in Serction 3 we have derived the direct theorem for the identification capacity with input cost (Theorem 3.1), whereas in Section 4 we have derived the direct theorem for the channel resolvability with input cost (Theorem 4.1).

In this section, we establish the converse theorems for both of the identification capacity with input cost and the channel resolvability with input cost. To this end, we need the following key lemma which reveals the intrinsic relation between the identification capacity $D(\mu, \lambda, \Gamma | \mathbf{W})$ and the channel resolvability $S(\delta, \Gamma | \mathbf{W})$. The implication of this lemma is that the direct theorem for the channel resolvability gives the converse theorem for the identification capacity, and, on the other hand, the direct theorem for the identification capacity gives the converse theorem for the channel resolvability. It should be noted here that the this lemma holds only under the *finiteness* assumption of input alphabet $\mathcal{X}$, but *not* in full generality, although this assumption is to be dispensed with in the next section.

**Lemma 5.1** Let $\mathbf{W} = \{W^n\}_{n=1}^{\infty}$ be any general channel with *finite* input alphabet $\mathcal{X}$. Then, for all $\delta \geq 0$, $\mu \geq 0$, $\lambda \geq 0$ such that $\delta < 1 - \mu - \lambda$, it holds that

$$D(\mu, \lambda, \Gamma | \mathbf{W}) \leq S(\delta, \Gamma | \mathbf{W}). \tag{5.1}$$

*Proof:*

The proof basically parallels that of Theorem 9 in Han and Verdú [4], except for that, here, we have to take account of the cost constraint. □

**Theorem 5.1** Let $\mathbf{W} = \{W^n\}_{n=1}^{\infty}$ be any general channel with *finite* input alphabet $\mathcal{X}$. Then, for all $\varepsilon \geq 0$, $\delta \geq 0$, $\mu \geq 0$, $\lambda \geq 0$ such that

$$\varepsilon \leq \mu, \quad \varepsilon \leq \lambda, \quad \delta < 1 - \mu - \lambda,$$

it holds that

$$\begin{aligned}
\sup_{\mathbf{X} \in \mathcal{S}_\Gamma} \underline{I}(\mathbf{X}; \mathbf{Y}) \leq C_s(\varepsilon, \Gamma | \mathbf{W}) &\leq D(\mu, \lambda, \Gamma | \mathbf{W}) \\
&\leq S(\delta, \Gamma | \mathbf{W}) \\
&\leq \sup_{\mathbf{X} \in \mathcal{S}_\Gamma} \overline{I}(\mathbf{X}; \mathbf{Y}). \tag{5.2}
\end{aligned}$$

*Proof:* It follows from (2.7), (2.9), (4.3) and (5.1). □

An immediate consequence of Theorem 5.1 and Theorem 2.3 is:

**Theorem 5.2** If a general channel $\mathbf{W}$ with *finite* input alphabet $\mathcal{X}$ satisfies the strong converse property with input cost constraint $\Gamma$, then, for all $\varepsilon \geq 0$, $\delta \geq 0$, $\mu \geq 0$, $\lambda \geq 0$ such that

$$\varepsilon < 1, \ \mu + \lambda < 1, \ \delta < 1, \tag{5.3}$$

it holds that

$$\begin{aligned}\sup_{\mathbf{X} \in \mathcal{S}_\Gamma} \underline{I}(\mathbf{X};\mathbf{Y}) = C_s(\varepsilon, \Gamma|\mathbf{W}) &= D(\mu, \lambda, \Gamma|\mathbf{W}) \\ &= S(\delta, \Gamma|\mathbf{W}) \\ &= \sup_{\mathbf{X} \in \mathcal{S}_\Gamma} \overline{I}(\mathbf{X};\mathbf{Y}).\end{aligned} \tag{5.4}$$

From Theorem 5.2 we have the following two corollaries.

**Corollary 5.1** If a general channel $\mathbf{W}$ with *finite* input alphabet $\mathcal{X}$ satisfies the strong converse property with input cost constraint $\Gamma$, then, for all $\mu \geq 0$, $\lambda \geq 0$ such that $\mu + \lambda < 1$, it holds that

$$D(\mu, \lambda, \Gamma|\mathbf{W}) = C_s(\Gamma|\mathbf{W}).$$

**Corollary 5.2** If a general channel $\mathbf{W}$ with *finite* input alphabet $\mathcal{X}$ satisfies the strong converse property with input cost constraint $\Gamma$, then, for all $\delta$ such that $0 \leq \delta < 1$, it holds that

$$S(\delta, \Gamma|\mathbf{W}) = C_s(\Gamma|\mathbf{W}).$$

**Example 5.1** Let us consider any stationary memoryless channel $\mathbf{W} = \{W^n\}_{n=1}^\infty$ with finite input and output alphabets $\mathcal{X}, \mathcal{Y}$, which can be specified simply by the transition probability matrix $W : \mathcal{X} \to \mathcal{Y}$. Moreover, let the cost function $\mathbf{c} = \{c_n\}_{n=1}^\infty$ be *additive*, i.e., suppose that there exists a function $c : \mathcal{X} \to \mathbf{R}$ such that $c_n(\mathbf{x}) = c(x_1) + c(x_2) + \cdots + c(x_n)$ for $\mathbf{x} = (x_1, x_2, \cdots, x_n) \in \mathcal{X}^n$. It is shown in Han [6] that such a channel $W : \mathcal{X} \to \mathcal{Y}$ satisfies the strong converse property with any input cost constraint $\Gamma$. Therefore, by means of Corollary 5.1 and Corollary 5.2 it is concluded that, for all $\mu \geq 0$, $\lambda \geq 0$ ($\mu + \lambda < 1$) and all $0 \leq \delta < 1$, it holds that

$$D(\mu, \lambda, \Gamma|\mathbf{W}) = S(\delta, \Gamma|\mathbf{W}) = C_s(\Gamma|\mathbf{W}) = \max_{X: \mathrm{E}c(X) \leq \Gamma} I(X;Y),$$

where $I(X;Y)$ is the mutual information (cf. Cover and Thomas [12]) with the conditional probability distribution of $Y$ given $X$ fixed to $W$. □

**Example 5.2** Let us consider any stationary finite-memory channel $\mathbf{W} = \{W^n\}_{n=1}^{\infty}$ with finite input and output alphabets $\mathcal{X}, \mathcal{Y}$ (cf. Feinstein [7], Wolfowitz [8]). Moreover, let the cost function $\mathbf{c} = \{c_n\}_{n=1}^{\infty}$ be *additive*. Since this channel satisfies the strong converse property with any input cost constraint $\Gamma$ (cf. Wolfowitz [8], Han [6]), by means of Corollary 5.1 and Corollary 5.2 it is concluded that, for all $\mu \geq 0$, $\lambda \geq 0$ ($\mu + \lambda < 1$) and all $0 \leq \delta < 1$, it holds that

$$\begin{aligned} D(\mu, \lambda, \Gamma|\mathbf{W}) &= S(\delta, \Gamma|\mathbf{W}) = C_s(\Gamma|\mathbf{W}) \\ &= \liminf_{n \to \infty} \max_{X^n: \frac{1}{n}\mathrm{E}c_n(X^n) \leq \Gamma} \frac{1}{n} I(X^n; Y^n). \end{aligned}$$

**Example 5.3** Let us consider any irreducible unifilar finite-state channel $\mathbf{W} = \{W^n\}_{n=1}^{\infty}$ with finite input and output alphabets $\mathcal{X}, \mathcal{Y}$ (cf. Wolfowitz [9]). Moreover, let the cost function $\mathbf{c} = \{c_n\}_{n=1}^{\infty}$ be *additive*. Since this channel satisfies the strong converse property with any input cost constraint $\Gamma$ (cf. Wolfowitz [9], Han [6]), by means of Corollary 5.1 and Corollary 5.2 it is concluded that, for all $\mu \geq 0$, $\lambda \geq 0$ ($\mu + \lambda < 1$) and all $0 \leq \delta < 1$, it holds that

$$\begin{aligned} D(\mu, \lambda, \Gamma|\mathbf{W}) &= S(\delta, \Gamma|\mathbf{W}) = C_s(\Gamma|\mathbf{W}) \\ &= \liminf_{n \to \infty} \max_{X^n: \frac{1}{n}\mathrm{E}c_n(X^n) \leq \Gamma} \frac{1}{n} I(X^n; Y^n|s_0), \end{aligned}$$

where $s_0$ is the fixed initial state. □

# 6 Identification capacity and resolvability for channels with continuous input alphabet

Thus far, we have established the fundamental theorem (Theorem 5.2) on the identification capacity as well as on the channel resolvability. Here, it shoulde be noted that the validity of Theorem 5.2 is based on Lemma 5.1 that heavily depends on the *finiteness* of input alphabet $\mathcal{X}$. As a result, the validity of Theorem 5.2 also heavily relies upon the *finiteness* of input alphabet $\mathcal{X}$.

However, under some mild conditions on the *stability* of the distribution of the channel output with respect to channel inputs, we can dispense with this finiteness assumption. In this section, as an illustrative case, we consider the general channel $\mathbf{W}$ with $\mathbf{R}$ (the set of real numbers) as the input alphabet $\mathcal{X}$, i.e., $\mathcal{X} = \mathbf{R}$. We first establish the fundamental theorem for this case, and after that, as its applications, we show the formulas for the stationary additive white (and also non-white) Gaussian noise channels.

Now, let $\mathbf{W} = \{W^n\}_{n=1}^{\infty}$ be any general channel with input alphabet $\mathcal{X} = \mathbf{R}$ and arbitrary output alphabet $\mathcal{Y}$. We assume the following two properties. First, as in (2.4), put

$$\mathcal{X}^n(\Gamma) = \left\{ \mathbf{x} \in \mathcal{X}^n \,\middle|\, \frac{1}{n} c_n(\mathbf{x}) \leq \Gamma \right\}, \tag{6.1}$$

and assume that there exists an $n$-dimensional cube $V_n(\Gamma) \subset \mathcal{X}^n$ of edge length $l_n(\Gamma) \geq 2$ such that

$$\limsup_{n \to \infty} \frac{1}{n} \log \log l_n(\Gamma) = 0 \tag{6.2}$$

and

$$\mathcal{X}^n(\Gamma) \subset V_n(\Gamma) \quad (\forall n = 1, 2, \cdots). \tag{6.3}$$

Next, let $D(W^n(\,\cdot\,|\mathbf{v})\|W^n(\,\cdot\,|\mathbf{x}))$ $(\mathbf{v}, \mathbf{x} \in \mathcal{X}^n)$ denote the divergence (cf. Csiszár and Körner [11]) between $W^n(\,\cdot\,|\mathbf{v})$ and $W^n(\,\cdot\,|\mathbf{x})$, and define the $n \times n$ Fisher information matrix $F_n(\mathbf{x})$ by

$$F_n(\mathbf{x}) = \left[ \frac{\partial^2}{\partial v_i \partial v_j} D(W^n(\,\cdot\,|\mathbf{v})\|W^n(\,\cdot\,|\mathbf{x})) \right]_{\mathbf{v}=\mathbf{x}}, \tag{6.4}$$

where $\mathbf{v} = (v_1, v_2, \cdots, v_n)$. Moreover, define the norm of $F_n(\mathbf{x})$ by

$$\|F_n(\mathbf{x})\| = \sup_{\mathbf{v} \in \mathcal{X}^n, \mathbf{v} \neq \mathbf{0}} \frac{\mathbf{v} F_n(\mathbf{x}) \mathbf{v}^{\mathrm{T}}}{\mathbf{v}\mathbf{v}^{\mathrm{T}}},$$

where "T" denote the transpose of a matrix. Here, we assume that

$$\limsup_{n \to \infty} \frac{1}{n} \log \max \left\{ 1, \log \left( \sup_{\mathbf{x} \in V_n(\Gamma)} \|F_n(\mathbf{x})\| \right) \right\} = 0. \tag{6.5}$$

Notice that assumptions (6.2), (6.5) are rather weaker ones.

Under these assumptions, we make the following operation. First, (6.5) can be rewritten as

$$\sup_{\mathbf{x} \in V_n(\Gamma)} ||F_n(\mathbf{x})|| \leq \exp(e^{n\delta_n}), \tag{6.6}$$

where $\{\delta_n \geq 0\}_{n=1}^{\infty}$ is some sequence such that $\lim_{n \to \infty} \delta_n = 0$. Set

$$\gamma_n = \max(\delta_n, \frac{1}{\sqrt{n}}), \tag{6.7}$$

and divide the cube $V_n(\Gamma)$ in (6.3) into $k_n(\Gamma)$ finer cubes $\Lambda_n^{(i)}$ ($i = 1, 2, \cdots, k_n(\Gamma)$) each of edge length

$$\Delta_n = \exp(-e^{n\gamma_n}), \tag{6.8}$$

where the number $k_n(\Gamma)$ of these finer cubes is given by

$$k_n(\Gamma) = \left(\frac{l_n(\Gamma)}{\Delta_n}\right)^n.$$

It then follows from (6.2), (6.8) and $\lim_{n \to \infty} \gamma_n = 0$ that

$$\limsup_{n \to \infty} \frac{1}{n} \log \log k_n(\Gamma) = 0. \tag{6.9}$$

For each $i = 1, 2, \cdots, k_n(\Gamma)$, fix a representative point $\mathbf{u}_i$ in the cube $\Lambda_n^{(i)}$, and put

$$\mathcal{R}_n(\Gamma) = \left\{\mathbf{u}_1, \mathbf{u}_2, \cdots, \mathbf{u}_{k_n(\Gamma)}\right\}. \tag{6.10}$$

Furthermore, to each probability distribution $Q$ on $\mathcal{X}^n$ such that $Q(\mathcal{X}^n(\Gamma)) = 1$ we uniquely correspond the probability distribution $\overline{Q}$ on $\mathcal{R}_n(\Gamma)$ so that

$$\overline{Q}(\mathbf{u}_i) = Q(\Lambda_n^{(i)}) \quad (i = 1, 2, \cdots, k_n(\Gamma)). \tag{6.11}$$

We call this operation the *quantization* of probability distributions. Then, we have the following lemma to evaluate the variational distance $d(QW^n, \overline{Q}W^n)$ between probability distributions $QW^n$ and $\overline{Q}W^n$.

**Lemma 6.1** Under assumptions (6.2), (6.5), for any probability distribution $Q$ such that $Q(\mathcal{X}^n(\Gamma)) = 1$ and for all $n = 1, 2, \cdots$ it holds that

$$d(QW^n, \overline{Q}W^n) \leq \sqrt{n} \exp(-\frac{1}{2}e^{\sqrt{n}}). \tag{6.12}$$

*Proof:*

For any $\mathbf{x} \in \Lambda_n^{(i)}$ we use the Taylor expansion in $\mathbf{x}$ of the divergence $D(W^n(\,\cdot\,|\mathbf{x})\|W^n(\,\cdot\,|\mathbf{u}_i))$ to obtain

$$D(W^n(\,\cdot\,|\mathbf{x})\|W^n(\,\cdot\,|\mathbf{u}_i))$$
$$= \frac{1}{2}(\mathbf{x} - \mathbf{u}_i) F_n(\mathbf{x}')(\mathbf{x} - \mathbf{u}_i)^{\mathrm{T}},$$

where $\mathbf{x}' = \mathbf{u}_i + \theta(\mathbf{x} - \mathbf{u}_i)$ ($0 < \exists \theta < 1$). Notice that $\mathbf{x}' \in \Lambda_n^{(i)}$. Then, by the definition of the norm $\|F_n(\mathbf{x}')\|$ we have

$$D(W^n(\,\cdot\,|\mathbf{x})\|W^n(\,\cdot\,|\mathbf{u}_i))$$
$$\leq \frac{1}{2} \|F_n(\mathbf{x}')\|(\mathbf{x} - \mathbf{u}_i)(\mathbf{x} - \mathbf{u}_i)^{\mathrm{T}}. \tag{6.13}$$

We observe here that, since $\mathbf{x}$ and $\mathbf{u}_i$ belong to $\Lambda_n^{(i)}$, the modulus of each component of the vector $\mathbf{x} - \mathbf{u}_i$ is less than or equal to $\Delta_n$. Hence,

$$(\mathbf{x} - \mathbf{u}_i)(\mathbf{x} - \mathbf{u}_i)^{\mathrm{T}} \leq n\Delta_n^2. \tag{6.14}$$

Therefore, by means of (6.6)$\sim$ (6.8) and (6.14), (6.13) yields

$$D(W^n(\,\cdot\,|\mathbf{x})\|W^n(\,\cdot\,|\mathbf{u}_i))$$
$$\leq \frac{n}{2} \exp(e^{n\delta_n}) \exp(-2e^{n\gamma_n})$$
$$= \frac{n}{2} \exp(e^{n\delta_n} - 2e^{n\gamma_n})$$
$$\leq \frac{n}{2} \exp(-e^{n\gamma_n})$$
$$\leq \frac{n}{2} \exp(-e^{\sqrt{n}}). \tag{6.15}$$

Moreover, if we use the inequality (cf. Csiszár and Körner [11]):

$$\frac{1}{2} d^2(P_1, P_2) \leq D(P_1\|P_2),$$

it follows from (6.15) that

$$d(W^n(\,\cdot\,|\mathbf{x}), W^n(\,\cdot\,|\mathbf{u}_i)) \leq \sqrt{n} \exp(-\frac{1}{2} e^{\sqrt{n}}).$$

As a consequence, for any $\mathbf{x} \in \Lambda_n^{(i)}$ and for any subset $\mathcal{B} \subset \mathcal{Y}^n$,

$$|W^n(\mathcal{B}|\mathbf{x}) - W^n(\mathcal{B}|\mathbf{u}_i)| \leq \frac{\sqrt{n}}{2} \exp(-\frac{1}{2} e^{\sqrt{n}}). \tag{6.16}$$

On the other hand, since $Q(\mathcal{X}^n(\Gamma)) = 1$, from (6.11) we have

$$QW^n(\mathcal{B}) - \overline{Q}W^n(\mathcal{B})$$
$$= \sum_{\mathbf{x} \in \mathcal{X}^n(\Gamma)} W^n(\mathcal{B}|\mathbf{x}) Q(\mathbf{x}) - \sum_{i=1}^{k_n(\Gamma)} W^n(\mathcal{B}|\mathbf{u}_i) \overline{Q}(\mathbf{u}_i)$$
$$= \sum_{i=1}^{k_n(\Gamma)} \sum_{\mathbf{x} \in \Lambda_n^{(i)}} W^n(\mathcal{B}|\mathbf{x}) Q(\mathbf{x}) - \sum_{i=1}^{k_n(\Gamma)} \sum_{\mathbf{x} \in \Lambda_n^{(i)}} W^n(\mathcal{B}|\mathbf{u}_i) Q(\mathbf{x})$$
$$= \sum_{i=1}^{k_n(\Gamma)} \sum_{\mathbf{x} \in \Lambda_n^{(i)}} (W^n(\mathcal{B}|\mathbf{x}) - W^n(\mathcal{B}|\mathbf{u}_i)) Q(\mathbf{x}).$$

Hence,

$$|QW^n(\mathcal{B}) - \overline{Q}W^n(\mathcal{B})|$$
$$\leq \sum_{i=1}^{k_n(\Gamma)} \sum_{\mathbf{x} \in \Lambda_n^{(i)}} |W^n(\mathcal{B}|\mathbf{x}) - W^n(\mathcal{B}|\mathbf{u}_i)| Q(\mathbf{x}). \qquad (6.17)$$

Substitution of (6.16) into the right-hand side of (6.17) gives

$$|QW^n(\mathcal{B}) - \overline{Q}W^n(\mathcal{B})| \leq \sum_{i=1}^{k_n(\Gamma)} \overline{Q}(\mathbf{u}_i) \frac{\sqrt{n}}{2} \exp(-\frac{1}{2} e^{\sqrt{n}})$$
$$= \frac{\sqrt{n}}{2} \exp(-\frac{1}{2} e^{\sqrt{n}}).$$

Thus, noting that $\mathcal{B} \subset \mathcal{Y}^n$ is arbitrary, we have

$$d(QW^n, \overline{Q}W^n) = 2 \sup_{\mathcal{B} \subset \mathcal{Y}^n} |QW^n(\mathcal{B}) - \overline{Q}W^n(\mathcal{B})|$$
$$\leq \sqrt{n} \exp(-\frac{1}{2} e^{\sqrt{n}}),$$

thereby proving the claim of the lemma. □

With these preparations, in the sequel let us show formulas for the identification capacity and the channel resolvability for a general channel $\mathbf{W}$ with continuous input alphabet $\mathcal{X} = \mathbf{R}$ under input cost constraint $\Gamma$. To this end, we first need the following key lemma, which is the continuous-input counterpart of Lemma 5.1.

**Lemma 6.2** Let **W** be any general channel with input alphabet $\mathcal{X} = \mathbf{R}$ that satisfies assumptions (6.2), (6.5). Then, for all $\delta \geq 0$, $\mu \geq 0$, $\lambda \geq 0$ such that $\delta < 1 - \mu - \lambda$, it holds that

$$D(\mu, \lambda, \Gamma | \mathbf{W}) \leq S(\delta, \Gamma | \mathbf{W}). \tag{6.18}$$

*Proof:*

Let $R_1$ be any $(\mu, \lambda, \Gamma)$-achievable identification code rate. Then, by the definition of the achievability, there exists an $(n, N_n, \mu_n, \lambda_n, \Gamma)$ code such that

$$\limsup_{n \to \infty} \mu_n \leq \mu, \quad \limsup_{n \to \infty} \lambda_n \leq \lambda, \tag{6.19}$$

$$\liminf_{n \to \infty} \frac{1}{n} \log \log N_n \geq R_1. \tag{6.20}$$

Let $Q_1, Q_2, \cdots, Q_{N_n}$ (probability distributions on $\mathcal{X}^n(\Gamma)$) be the codewords of this identification code, and let $\mathcal{D}_1, \mathcal{D}_2, \cdots, \mathcal{D}_{N_n} \subset \mathcal{Y}^n$ be the corresponding decoding sets. Then, for any $j \neq k$,

$$\begin{aligned} d(Q_j W^n, Q_k W^n) &\geq 2(Q_j W^n(\mathcal{D}_j) - Q_k W^n(\mathcal{D}_j)) \\ &\geq 2(1 - \mu_n - \lambda_n). \end{aligned}$$

It then follows from (6.19) that

$$\begin{aligned} \liminf_{n \to \infty} d(Q_j W^n, Q_k W^n) &\geq 2(1 - \limsup_{n \to \infty} \mu_n - \limsup_{n \to \infty} \lambda_n) \\ &\geq 2(1 - \mu - \lambda). \end{aligned} \tag{6.21}$$

On the other hand, let $R_2$ be any $(\delta, \Gamma)$-achievable resolvability rate. Then, for all $j = 1, 2, \cdots, N_n$, there exists a deterministic encoder

$$\varphi_n^{(j)} : \mathcal{M}_n \equiv \{1, 2, \cdots, M_n\} \to \mathcal{X}^n(\Gamma)$$

such that

$$\limsup_{n \to \infty} d(Q_j W^n, \tilde{Q}_j W^n) \leq \delta, \tag{6.22}$$

$$\limsup_{n \to \infty} \frac{1}{n} \log M_n \leq R_2, \tag{6.23}$$

where $\tilde{Q}_j$ is the probability distribution (on $\mathcal{X}^n(\Gamma)$) of $\tilde{X}_j^n \equiv \varphi_n^{(j)}(U_{M_n})$ and $U_{M_n}$ is the uniform random number as specified in Section 4.

Now, let $\overline{Q}_j$ denote the quantized probability distribution $\overline{Q}$ on $\mathcal{R}_n(\Gamma)$ as specified by (6.11) with $\tilde{Q}_j$ instead of $Q$. We observe here from the way of constructing $\overline{Q}_j$ that there must exist a deterministic encoder $\overline{\phi}_n^{(j)} : \mathcal{M}_n \to \mathcal{R}_n(\Gamma)$ such that $\overline{Q}_j$ is the probability distribution (on $\mathcal{R}_n(\Gamma)$) of $\overline{X}_j^n \equiv \overline{\phi}_n^{(j)}(U_{M_n})$. Then, by virtue of Lemma 6.1 with $\tilde{Q}_j, \overline{Q}_j$ instead of $Q, \overline{Q}$ we have

$$d(\tilde{Q}_j W^n, \overline{Q}_j W^n) \leq \sqrt{n} \exp(-\frac{n}{2} e^{\sqrt{n}}),$$

from which follows that

$$\lim_{n \to \infty} d(\tilde{Q}_j W^n, \overline{Q}_j W^n) = 0. \tag{6.24}$$

On the other hand, in view of

$$\begin{aligned} &d(Q_j W^n, \overline{Q}_j W^n) \\ &\leq d(Q_j W^n, \tilde{Q}_j W^n) + d(\tilde{Q}_j W^n, \overline{Q}_j W^n), \end{aligned}$$

(6.22) and (6.24) yield

$$\limsup_{n \to \infty} d(Q_j W^n, \overline{Q}_j W^n) \leq \delta \quad (j = 1, 2, \cdots, N_n). \tag{6.25}$$

To argue otherwise, suppose that $\overline{Q}_j = \overline{Q}_k$ for some $j \neq k$. Then, since

$$\begin{aligned} d(Q_j W^n, Q_k W^n) &\leq d(Q_j W^n, \overline{Q}_j W^n) + d(Q_k W^n, \overline{Q}_j W^n) \\ &= d(Q_j W^n, \overline{Q}_j W^n) + d(Q_k W^n, \overline{Q}_k W^n), \end{aligned}$$

(6.25) gives

$$\limsup_{n \to \infty} d(Q_j W^n, Q_k W^n) \leq 2\delta. \tag{6.26}$$

However, (6.26) contradicts (6.21), because of the assumption $\delta < 1 - \mu - \lambda$. Thus, if $j \neq k$, then $\overline{Q}_j \neq \overline{Q}_k$ must hold. Therefore, all the distributions on $\mathcal{R}_n(\Gamma)$:

$$\overline{Q}_1, \overline{Q}_2, \cdots, \overline{Q}_{N_n}$$

must be different from one another. We notice that the distribution $\overline{Q}_j$ on $\mathcal{R}_n(\Gamma)$ can be expressed in the form as

$$\overline{Q}_j(\mathbf{u}) = \frac{m}{M_n} \quad (\forall \mathbf{u} \in \mathcal{R}_n(\Gamma); \forall j = 1, 2, \cdots, N_n), \tag{6.27}$$

where $m$ is some nonnegative integer. Since the number of different probability distributions $\overline{Q}_j$ satisfying property (6.27) is at most

$$|\mathcal{R}_n(\Gamma)|^{M_n} = (k_n(\Gamma))^{M_n},$$

it must hold that $N_n \leq (k_n(\Gamma))^{M_n}$. Hence,

$$\log N_n \leq M_n \log k_n(\Gamma),$$

and so

$$\frac{1}{n} \log \log N_n \leq \frac{1}{n} \log M_n + \frac{1}{n} \log \log k_n(\Gamma).$$

Then, from (6.9) we have

$$\liminf_{n \to \infty} \frac{1}{n} \log \log N_n \leq \limsup_{n \to \infty} \frac{1}{n} \log M_n,$$

which, together with (6.20) and (6.23), implies that

$$R_1 \leq R_2.$$

Thus, noting that $R_1$ is an arbitrary $(\mu, \lambda, \Gamma)$-achievable identification code rate, and $R_2$ is an arbitrary $(\delta, \Gamma)$-achievable resolvability rate, we complete the proof of the lemma. □

Now, as the continuous-input counterpart of Theorem 5.1, we have the following theorem.

**Theorem 6.1** Let $\mathbf{W}$ be any general channel with input alphabet $\mathcal{X} = \mathbf{R}$ that satisfies assumptions (6.2), (6.5). Then, for all $\varepsilon \geq 0$, $\delta \geq 0$, $\mu \geq 0$, $\lambda \geq 0$ such that

$$\varepsilon \leq \mu, \quad \varepsilon \leq \lambda, \quad \delta < 1 - \mu - \lambda,$$

it holds that

$$\sup_{\mathbf{X} \in \mathcal{S}_\Gamma} \underline{I}(\mathbf{X}; \mathbf{Y}) \leq C_s(\varepsilon, \Gamma | \mathbf{W}) \leq D(\mu, \lambda, \Gamma | \mathbf{W})$$
$$\leq S(\delta, \Gamma | \mathbf{W})$$
$$\leq \sup_{\mathbf{X} \in \mathcal{S}_\Gamma} \overline{I}(\mathbf{X}; \mathbf{Y}). \qquad (6.28)$$

*Proof:* It follows from (2.7), (2.9), (4.3) and (6.18). □

An immediate consequence of Theorem 6.1 and Theorem 2.3 is:

**Theorem 6.2 (Fundamental theorem)** Let $\mathbf{W}$ be any general channel with input alphabet $\mathcal{X} = \mathbf{R}$ that satisfies assumptions (6.2), (6.5) as well as the strong converse property with input cost constraint $\Gamma$. Then, for all $\varepsilon \geq 0$, $\delta \geq 0$, $\mu \geq 0$, $\lambda \geq 0$ such that

$$\varepsilon < 1, \ \mu + \lambda < 1, \ \delta < 1,$$

it holds that

$$\sup_{\mathbf{X} \in \mathcal{S}_\Gamma} \underline{I}(\mathbf{X}; \mathbf{Y}) = C_s(\varepsilon, \Gamma | \mathbf{W}) = D(\mu, \lambda, \Gamma | \mathbf{W})$$
$$= S(\delta, \Gamma | \mathbf{W})$$
$$= \sup_{\mathbf{X} \in \mathcal{S}_\Gamma} \overline{I}(\mathbf{X}; \mathbf{Y}). \quad (6.29)$$

From Theorem 6.2 we have the following two corollaries.

**Corollary 6.1** Let $\mathbf{W}$ be any general channel with input alphabet $\mathcal{X} = \mathbf{R}$ that satisfies assumptions (6.2), (6.5) as well as the strong converse property with input cost constraint $\Gamma$. Then, for all $\mu \geq 0$, $\lambda \geq 0$ such that $\mu + \lambda < 1$, it holds that
$$D(\mu, \lambda, \Gamma | \mathbf{W}) = C_s(\Gamma | \mathbf{W}).$$

**Corollary 6.2** Let $\mathbf{W}$ be any general channel with input alphabet $\mathcal{X} = \mathbf{R}$ that satisfies assumptions (6.2), (6.5) as well as the strong converse property with input cost constraint $\Gamma$. Then, for all $\delta$ such that $0 \leq \delta < 1$, it holds that
$$S(\delta, \Gamma | \mathbf{W}) = C_s(\Gamma | \mathbf{W}).$$

## 6.1 AWGN channels

Now, as an illustrative typical application of Theorem 6.2 as well as Corollaries 6.1 and 6.2, let us first consider the stationary additive white Gaussian noise channel $\mathbf{W} = \{W : \mathcal{X} \to \mathcal{Y}\}$ (called the AWGN channel) with input and output alphabets $\mathcal{X} = \mathcal{Y} = \mathbf{R}$. The transition probability density of the AWGN channel with noise power $N > 0$ is given by

$$W(y|x) = \frac{1}{\sqrt{2\pi N}} e^{-\frac{(y-x)^2}{2N}}. \quad (6.30)$$

Set

$$\mathbf{v} = (v_1, v_2, \cdots, v_n) \in \mathcal{X}^n,$$
$$\mathbf{x} = (x_1, x_2, \cdots, x_n) \in \mathcal{X}^n,$$

then, by the memorylessness of the channel we have

$$D(W^n(\cdot|\mathbf{v})||W^n(\cdot|\mathbf{x}))$$
$$= \sum_{i=1}^{n} D(W(\cdot|v_i)||W(\cdot|x_i)).$$

It is easy to check that

$$D(W(\cdot|v_i)||W(\cdot|x_i)) = \frac{(v_i - x_i)^2}{2N}.$$

Hence,

$$D(W^n(\cdot|\mathbf{v})||W^n(\cdot|\mathbf{x})) = \sum_{i=1}^{n} \frac{(v_i - x_i)^2}{2N}.$$

Then, the Fisher information matrix $F_n(\mathbf{x})$ in (6.4) reduces to the diagonal matrix as

$$F_n(\mathbf{x}) = \begin{pmatrix} \frac{1}{N} & & & O \\ & \frac{1}{N} & & \\ & & \ddots & \\ O & & & \frac{1}{N} \end{pmatrix}.$$

Hence, we have $||F_n(\mathbf{x})|| = \frac{1}{N}$, which obviously satisfies assumption (6.5). Next, the input power constraint for the AWGN channel can be written as

$$\frac{1}{n}c_n(\mathbf{x}) \equiv \frac{1}{n}(x_1^2 + x_2^2 + \cdots x_n^2) \le P \quad (\forall n = 1, 2, \cdots),$$

where $\mathbf{x} = (x_1, x_2, \cdots, x_n) \in \mathcal{X}^n$ and $P > 0$ is the signal power. Therefore,

$$-\sqrt{nP} \le x_i \le \sqrt{nP} \quad (i = 1, 2, \cdots, n)$$

must hold. Hence,

$$\mathcal{X}^n(P) = \left\{ \mathbf{x} \in \mathcal{X}^n \,\bigg|\, \frac{1}{n}c_n(\mathbf{x}) \le P \right\}$$

is contained in an $n$-dimensional cube $V_n(P)$ of edge length $l_n(P) = 2\sqrt{nP}$. Then, it is obvious that assumption (6.2) is satisfied. Thus, it is concluded that both of assumptions (6.2) and (6.5) are satisfied for the AWGN channel. On the other hand, it is shown in Shannon [10] (also, cf. Han [6]) that the AWGN channel satisfies the strong converse property with any input cost constraint $P$. Thus, by means of Theorem 6.2 we have the following theorem.

**Theorem 6.3** Let $\mathbf{W}$ be the AWGN channel with signal power $P$ and noise power $N$. Then, for all $\varepsilon \geq 0$, $\delta \geq 0$, $\mu \geq 0$, $\lambda \geq 0$ such that

$$\varepsilon < 1, \ \mu + \lambda < 1, \ \delta < 1,$$

it holds that

$$\begin{aligned}
C_s(\varepsilon, P|\mathbf{W}) &= D(\mu, \lambda, P|\mathbf{W}) \\
&= S(\delta, P|\mathbf{W}) \\
&= C_s(P|\mathbf{W}) \\
&= \frac{1}{2}\log\left(1 + \frac{P}{N}\right).
\end{aligned} \quad (6.31)$$

**Corollary 6.3** Let $\mathbf{W}$ be the AWGN channel with signal power $P$ and noise power $N$. Then, for all $\mu \geq 0$, $\lambda \geq 0$ such that $\mu + \lambda < 1$, it holds that

$$D(\mu, \lambda, P|\mathbf{W}) = C_s(P|\mathbf{W}) = \frac{1}{2}\log\left(1 + \frac{P}{N}\right).$$

**Corollary 6.4** Let $\mathbf{W}$ be the AWGN channel with signal power $P$ and noise power $N$. Then, for all $\delta$ such that $0 \leq \delta < 1$, it holds that

$$S(\delta, P|\mathbf{W}) = C_s(P|\mathbf{W}) = \frac{1}{2}\log\left(1 + \frac{P}{N}\right).$$

## 6.2 ANWGN channels

We may consider a more general Gaussian channel, i.e., the stationary additive but *non*-white Gaussian noise channel (called the ANWGN channel) with input and output alphabets $\mathcal{X} = \mathcal{Y} = \mathbf{R}$. For this ANWGN channel too, the counterpart of Theorem 6.3 as well as Corollaries 6.3, 6.4 holds. In this subsection we show this in several steps. Let $\mathbf{Z} \equiv (Z_1, Z_2, \cdots)$ be the stationary non-white Gaussian noise process with mean zero. Then, the ANWGN channel is specified by

$$Y_i^{(n)} = X_i^{(n)} + Z_i \quad (i = 1, 2, \cdots, n), \quad (6.32)$$

where $X_i^{(n)}$ is the $i$-th channel input and $Y_i^{(n)}$ is the corresponding channel output, and $(Z_1, Z_2, \cdots, Z_n)$ is independent from $(X_1^{(n)}, X_2^{(n)}, \cdots, X_n^{(n)})$ (*additivity*). Throughout in this subsection we assume that the noise process $\mathbf{Z} \equiv (Z_1, Z_2, \cdots)$ is *purely nondeterministic* (cf. Ihara [14]). Let this channel be denoted by $\mathbf{W} = \{W^n : \mathcal{X}^n \to \mathcal{Y}^n\}_{n=1}^{\infty}$. Define the autocorrelations $\{\gamma_k\}_{k=-\infty}^{k=\infty}$ of the noise process $\mathbf{Z}$ by

$$\gamma_k = \mathrm{E}(Z_i Z_{i+k}) \quad (k = 0, 1, 2, \cdots), \tag{6.33}$$

$$\gamma_k \equiv \gamma_{-k}, \quad (k = -1, -2, \cdots), \tag{6.34}$$

where it should be noted that the right-hand side of (6.33) does not depend on $i$ because of the stationarity of the noise process. In terms of these $\gamma_k$'s we define the $n$-dimensional covariance matrix $V_n$ as

$$V_n = \begin{pmatrix} \gamma_0 & \gamma_1 & \gamma_2 & \cdots & \gamma_{n-1} \\ \gamma_{-1} & \gamma_0 & \gamma_1 & \cdots & \gamma_{n-2} \\ \gamma_{-2} & \gamma_{-1} & \gamma_0 & \cdots & \gamma_{n-3} \\ \cdots & \cdots & \ddots & & \\ \cdots & \cdots & & \ddots & \\ \gamma_{-(n-1)} & \gamma_{-(n-2)} & \gamma_{-(n-3)} & \cdots & \gamma_0 \end{pmatrix}. \tag{6.35}$$

With this $V_n$, the transition probability density of the channel $W^n$ is given by

$$W^n(\mathbf{y}|\mathbf{x}) = \frac{1}{\sqrt{(2\pi)^n \det V_n}} \exp\left[-\frac{1}{2}(\mathbf{y} - \mathbf{x}) V_n^{-1} (\mathbf{y} - \mathbf{x})^{\mathrm{T}}\right], \tag{6.36}$$

where $\mathbf{x} \in \mathcal{X}^n, \mathbf{y} \in \mathcal{Y}^n$ are an input and the corresponding output, respectively.

Let us now consider to equivalently transform the channel as specified by (6.32) as follows. We first note that, since $V_n$ is a symmetric positive-definite matrix, there exists an $n$-dimensional orthogonal martix $U_n$ such that

$$U_n^{\mathrm{T}} V_n U_n = \begin{pmatrix} N_1^{(n)} & & & O \\ & N_2^{(n)} & & \\ & & \ddots & \\ O & & & N_n^{(n)} \end{pmatrix}. \tag{6.37}$$

where $N_i^{(n)}$ ($i = 1, 2, \cdots, n$) are the eigen values of $V_n$. We notice here that $N_i^{(n)} > 0$ ($i = 1, 2, \cdots, n$) because the noise process $\mathbf{Z}$ is purely nondeterministic. Define the modified noise process $(\overline{Z}_1^{(n)}, \overline{Z}_2^{(n)}, \cdots, \overline{Z}_n^{(n)})$ by

$$(\overline{Z}_1^{(n)}, \overline{Z}_2^{(n)}, \cdots, \overline{Z}_n^{(n)}) = (Z_1, Z_2, \cdots, Z_n) U_n, \tag{6.38}$$

where it is evident that $\overline{Z}_1^{(n)}, \overline{Z}_2^{(n)}, \cdots, \overline{Z}_n^{(n)}$ are Gaussian and mutually independent with means zero and variances $\mathrm{E}(\overline{Z}_i^{(n)})^2 = N_i^{(n)}$ ($i = 1, 2, \cdots, n$). Accordingly, define the random variables $(\overline{X}_1^{(n)}, \overline{X}_2^{(n)}, \cdots, \overline{X}_n^{(n)})$ and $(\overline{Y}_1^{(n)}, \overline{Y}_2^{(n)}, \cdots, \overline{Y}_n^{(n)})$ by

$$(\overline{X}_1^{(n)}, \overline{X}_2^{(n)}, \cdots, \overline{X}_n^{(n)}) = (X_1^{(n)}, X_2^{(n)}, \cdots, X_n^{(n)})U_n, \qquad (6.39)$$

$$(\overline{Y}_1^{(n)}, \overline{Y}_2^{(n)}, \cdots, \overline{Y}_n^{(n)}) = (Y_1^{(n)}, Y_2^{(n)}, \cdots, Y_n^{(n)})U_n. \qquad (6.40)$$

Then, the channel as specified by (6.32) is equivalently transformed to the nonstationary but *memoryless* additive Gaussian channel specified by

$$\overline{Y}_i^{(n)} = \overline{X}_i^{(n)} + \overline{Z}_i^{(n)} \quad (i = 1, 2, \cdots, n). \qquad (6.41)$$

Let this channel be denoted by $\overline{\mathbf{W}} = \{\overline{W}^n : \mathcal{X}^n \to \mathcal{Y}^n\}_{n=1}^\infty$. We notice here that under this transformation we have

$$(\overline{X}_1^{(n)})^2 + (\overline{X}_2^{(n)})^2 + \cdots + (\overline{X}_n^{(n)})^2$$
$$= (X_1^{(n)})^2 + (X_2^{(n)})^2 + \cdots + (X_n^{(n)})^2. \qquad (6.42)$$

Therefore, setting

$$c_n(\mathbf{x}) = x_1^2 + x_2^2 + \cdots + x_n^2 \quad (\mathbf{x} = (x_1, x_2, \cdots, x_n) \in \mathcal{X}^n),$$

$$X^n = (X_1^{(n)}, X_2^{(n)}, \cdots, X_n^{(n)}),$$

$$\overline{X}^n = (\overline{X}_1^{(n)}, \overline{X}_2^{(n)}, \cdots, \overline{X}_n^{(n)}),$$

and noting that the transformation (6.39), (6.40) preserves the mutual information, we have, for $P > 0$,

$$\frac{1}{n} \max_{X^n : \frac{1}{n}\mathrm{E}c_n(X^n) \leq P} I(X^n; Y^n) = \frac{1}{n} \max_{\overline{X}^n : \frac{1}{n}\mathrm{E}c_n(\overline{X}^n) \leq P} I(\overline{X}^n; \overline{Y}^n), \qquad (6.43)$$

where $Y^n, \overline{Y}^n$ are the channel outputs via channels $W^n, \overline{W}^n$ due to the channel inputs $X^n, \overline{X}^n$, respectively. On the other hand, it is well-known that the right-hand side of (6.43) is explicitly written by using the technique of water-filling (e.g., see Gallager [13]) as

$$\frac{1}{n} \max_{\overline{X}^n : \frac{1}{n}\mathrm{E}c_n(\overline{X}^n) \leq P} I(\overline{X}^n; \overline{Y}^n) = \frac{1}{2n} \sum_{i=1}^n \log\left(1 + \frac{P_i^{(n)}}{N_i^{(n)}}\right), \qquad (6.44)$$

where
$$P_i^{(n)} = \max[A_P^{(n)} - N_i^{(n)}, 0] \quad (i = 1, 2, \cdots, n), \tag{6.45}$$

and $A_P^{(n)} > 0$ is specified by the equation

$$\sum_{i=1}^{n} P_i^{(n)} = nP. \tag{6.46}$$

Here, it is easy to see from (6.45) and (6.46) that

$$A_P^{(n)} \geq P. \tag{6.47}$$

Now, we have the following lemma.

**Lemma 6.3** (Ihara[14]) The $P$-cost capacity $C_s(P|\mathbf{W})$ of the ANWGN channel $\mathbf{W}$ with input cost constraiint

$$\frac{1}{n}(x_1^2 + x_2^2 + \cdots + x_n^2) \leq P \quad ((x_1, x_2, \cdots, x_n) \in \mathcal{X}^n) \tag{6.48}$$

is given by

$$\begin{aligned}
C_s(P|\mathbf{W}) &= \lim_{n \to \infty} \frac{1}{n} \max_{X^n : \frac{1}{n}\mathrm{E}c_n(X^n) \leq P} I(X^n; Y^n) \\
&= \lim_{n \to \infty} \frac{1}{2n} \sum_{i=1}^{n} \log\left(1 + \frac{P_i^{(n)}}{N_i^{(n)}}\right).
\end{aligned} \tag{6.49}$$

It is possible also to give the non-limiting formula for the right-hand side of (6.49). To do so, let us define the spectral density function $g(\lambda)$ of the noise process $\mathbf{Z} \equiv (Z_1, Z_2, \cdots)$ by

$$g(\lambda) = \frac{1}{2\pi} \sum_{k=-\infty}^{\infty} \gamma_k e^{-ik\lambda} \quad (-\pi \leq \lambda \leq \pi), \tag{6.50}$$

where $\gamma_k$ are the autocorrelations as defined in (6.33) and (6.34). Then, Lemma 6.3 can be rewritten as

**Lemma 6.4** (Ihara[14]) The $P$-cost capacity $C_s(P|\mathbf{W})$ of the ANWGN channel $\mathbf{W}$ with input cost constraiint (6.48) is given by

$$C_s(P|\mathbf{W}) = \frac{1}{4\pi} \int_{-\pi}^{\pi} \log\left(1 + \frac{f(\lambda)}{g(\lambda)}\right) d\lambda, \tag{6.51}$$

where
$$f(\lambda) = \max[\alpha_P - g(\lambda), 0] \quad (-\pi \le \lambda \le \pi), \qquad (6.52)$$
and $\alpha_P > 0$ is specified by the equation
$$\int_{-\pi}^{\pi} f(\lambda) d\lambda = P. \qquad (6.53)$$

We need also the following lemma.

**Lemma 6.5** The ANWGN channel $\mathbf{W}$ with input cost constraint (6.48) satisfies the strong converse property.

*Proof:*
In view of (6.42) and the equivalence of the channels $\mathbf{W}$ and $\overline{\mathbf{W}}$, it suffices to show the strong converse property of the channel $\overline{\mathbf{W}}$ under input cost constraint (6.48). To this end, in the light of Theorem 2.2 and Theorem 2.3 with $P$ instead of $\Gamma$, it suffices to show
$$\sup_{\overline{\mathbf{X}} \in \mathcal{S}_P} \overline{I}(\overline{\mathbf{X}}; \overline{\mathbf{Y}}) \le C_s(P|\mathbf{W}),$$
or equivalently (cf. Lemma 6.3),
$$\sup_{\overline{\mathbf{X}} \in \mathcal{S}_P} \overline{I}(\overline{\mathbf{X}}; \overline{\mathbf{Y}}) \le \lim_{n \to \infty} \frac{1}{2n} \sum_{i=1}^{n} \log\left(1 + \frac{P_i^{(n)}}{N_i^{(n)}}\right), \qquad (6.54)$$
where $\overline{\mathbf{X}}$ denotes an input process for the channel $\overline{\mathbf{W}}$ and $\overline{\mathbf{Y}}$ denotes the channel output process via $\overline{\mathbf{W}}$ due to $\overline{\mathbf{X}}$. First, let $\overline{\mathbf{X}} = \{\overline{X}^n\}_{n=1}^{\infty}$ be an arbitrary input such that $\overline{\mathbf{X}} \in \mathcal{S}_P$ and $\overline{\mathbf{Y}} = \{\overline{Y}^n\}_{n=1}^{\infty}$ be the corresponding output via $\overline{\mathbf{W}} = \{\overline{W}^n\}_{n=1}^{\infty}$ due to the input $\overline{\mathbf{X}}$. For simplicity, set
$$i(\overline{X}^n; \overline{Y}^n) = \frac{1}{n} \log \frac{\overline{W}^n(\overline{Y}^n|\overline{X}^n)}{P_{\overline{Y}^n}(\overline{Y}^n)}, \qquad (6.55)$$
which we transform as
$$i(\overline{X}^n; \overline{Y}^n) = \frac{1}{n} \log \frac{\overline{W}^n(\overline{Y}^n|\overline{X}^n)}{P_{\tilde{Y}^n}(\overline{Y}^n)} - \frac{1}{n} \log \frac{P_{\overline{Y}^n}(\overline{Y}^n)}{P_{\tilde{Y}^n}(\overline{Y}^n)}, \qquad (6.56)$$
where
$$P_{\tilde{Y}^n}(\mathbf{y}) = \prod_{i=1}^{n} P_{\tilde{Y}_i}(y_i) \quad (\mathbf{y} = (y_1, y_2, \cdots, y_n) \in \mathcal{Y}^n), \qquad (6.57)$$

and $P_{\tilde{Y}_i}$ is the probability density of the output $\overline{Y}_i^{(n)}$ via the $i$-th component Gaussian channel in (6.41) due to the input $\overline{X}_i^{(n)}$ that attains the maximum of the mutual information $I(\overline{X}_i^{(n)}; \overline{Y}_i^{(n)})$ under the condition $\mathrm{E}(\overline{X}_i^{(n)})^2 \leq P_i^{(n)}$ ($i = 1, 2, \cdots, n$). Specifically, for $i = 1, 2, \cdots, n$,

$$P_{\tilde{Y}_i}(y) = \frac{1}{\sqrt{2\pi(P_i^{(n)} + N_i^{(n)})}} e^{-\frac{y^2}{2(P_i^{(n)} + N_i^{(n)})}} \quad (y \in \mathcal{Y}). \tag{6.58}$$

It is easy to check that

$$\text{p-}\liminf_{n \to \infty} \frac{1}{n} \log \frac{P_{\overline{Y}^n}(\overline{Y}^n)}{P_{\tilde{Y}^n}(\overline{Y}^n)} \geq 0,$$

and hence, from (6.56),

$$\begin{aligned}
\text{p-}\limsup_{n \to \infty} i(\overline{X}^n; \overline{Y}^n) &\leq \text{p-}\limsup_{n \to \infty} \frac{1}{n} \log \frac{\overline{W}^n(\overline{Y}^n | \overline{X}^n)}{P_{\tilde{Y}^n}(\overline{Y}^n)} \\
&\quad - \text{p-}\liminf_{n \to \infty} \frac{1}{n} \log \frac{P_{\overline{Y}^n}(\overline{Y}^n)}{P_{\tilde{Y}^n}(\overline{Y}^n)} \\
&\leq \text{p-}\limsup_{n \to \infty} \frac{1}{n} \log \frac{\overline{W}^n(\overline{Y}^n | \overline{X}^n)}{P_{\tilde{Y}^n}(\overline{Y}^n)}. \tag{6.59}
\end{aligned}$$

On the other hand, putting

$$\begin{aligned}
\overline{X}^n &= (\overline{X}_1^{(n)}, \overline{X}_2^{(n)}, \cdots, \overline{X}_n^{(n)}), \\
\overline{Y}^n &= (\overline{Y}_1^{(n)}, \overline{Y}_2^{(n)}, \cdots, \overline{Y}_n^{(n)})
\end{aligned}$$

and noting (6.57) as well as the memorylessness of the channel $\overline{W}^n$, we have

$$\frac{1}{n} \log \frac{\overline{W}^n(\overline{Y}^n | \overline{X}^n)}{P_{\tilde{Y}^n}(\overline{Y}^n)} = \frac{1}{n} \sum_{i=1}^n \log \frac{\overline{W}_i(\overline{Y}_i^{(n)} | \overline{X}_i^{(n)})}{P_{\tilde{Y}_i}(\overline{Y}_i^{(n)})}, \tag{6.60}$$

where, in view of (6.41), for $i = 1, 2, \cdots, n$,

$$\overline{W}_i(y|x) = \frac{1}{\sqrt{2\pi N_i^{(n)}}} e^{-\frac{(y-x)^2}{2N_i^{(n)}}} \quad (x \in \mathcal{X}, y \in \mathcal{Y}). \tag{6.61}$$

Now, fix any realization $\mathbf{x} = (x_1, x_2, \cdots, x_n)$ of $\overline{X}^n$ and set

$$I(\overline{Y}_i^{(n)}|x_i) \equiv \log \frac{\overline{W}_i(\overline{Y}_i^{(n)}|x_i)}{P_{\tilde{Y}_i}(\overline{Y}_i^{(n)})}, \qquad (6.62)$$

then, $I(\overline{Y}_i^{(n)}|x_i)$ ($i = 1, 2, \cdots, n$) are mutually independent under the conditional distribution $\overline{W}^n(\cdot|\mathbf{x})$ given $\overline{X}^n = \mathbf{x}$, owing to the memorylessness of the channel $\overline{W}^n$. Then, it follows from (6.58) and (6.61) that

$$\log \frac{\overline{W}_i(y|x)}{P_{\tilde{Y}_i}(y)} = \frac{1}{2} \log\left(1 + \frac{P_i^{(n)}}{N_i^{(n)}}\right) + \frac{y^2}{2(P_i^{(n)} + N_i^{(n)})} - \frac{(y-x)^2}{2N_i^{(n)}}.$$

Therefore,

$$\mathrm{E}_{\mathbf{x}}\left[\log \frac{\overline{W}_i(\overline{Y}_i^{(n)}|x_i)}{P_{\tilde{Y}_i}(\overline{Y}_i^{(n)})}\right]$$
$$= \int_{-\infty}^{\infty} \overline{W}_i(y|x_i) \log \frac{\overline{W}_i(y|x_i)}{P_{\tilde{Y}_i}(y)} dy$$
$$= \frac{1}{2} \log\left(1 + \frac{P_i^{(n)}}{N_i^{(n)}}\right) + \frac{x_i^2 - P_i^{(n)}}{2(P_i^{(n)} + N_i^{(n)})},$$

where $\mathrm{E}_{\mathbf{x}}$ denotes the conditional expectation under the conditional distribution $\overline{W}^n(\cdot|\mathbf{x})$. Hence,

$$\mathrm{E}_{\mathbf{x}}\left[\frac{1}{n}\sum_{i=1}^n \log \frac{\overline{W}_i(\overline{Y}_i^{(n)}|x_i)}{P_{\tilde{Y}_i}(\overline{Y}_i^{(n)})}\right]$$
$$= \frac{1}{2n} \sum_{i=1}^n \log\left(1 + \frac{P_i^{(n)}}{N_i^{(n)}}\right) + \sum_{i=1}^n \frac{x_i^2 - P_i^{(n)}}{2n(P_i^{(n)} + N_i^{(n)})}$$
$$\leq \frac{1}{2n} \sum_{i=1}^n \log\left(1 + \frac{P_i^{(n)}}{N_i^{(n)}}\right) + \frac{\sum_{i=1}^n x_i^2 - nP}{2nA_P^{(n)}},$$

where $A_P^{(n)} > 0$ is as specified by (6.45), (6.46), and we have taken account of the fact that $P_i^{(n)} > 0$ (i.e., $A_P^{(n)} > N_i^{(n)}$) implies $P_i^{(n)} + N_i^{(n)} = A_P^{(n)}$ and $P_i^{(n)} = 0$ implies $N_i^{(n)} \geq A_P^{(n)}$. Then, in view of input cost constraint (6.48), we have

$$\mathrm{E}_{\mathbf{x}}\left[\frac{1}{n}\sum_{i=1}^n \log \frac{\overline{W}_i(\overline{Y}_i^{(n)}|x_i)}{P_{\tilde{Y}_i}(\overline{Y}_i^{(n)})}\right] \leq \frac{1}{2n} \sum_{i=1}^n \log\left(1 + \frac{P_i^{(n)}}{N_i^{(n)}}\right). \qquad (6.63)$$

On the other hand, a direct calculation shows that the variance of $I(\overline{Y}_i^{(n)}|x_i)$ under the conditional distribution $\overline{W}^n(\cdot|\mathbf{x})$ is given by

$$V_{\mathbf{x}}\left[\log \frac{\overline{W}_i(\overline{Y}_i^{(n)}|x_i)}{P_{\tilde{Y}_i}(\overline{Y}_i^{(n)})}\right]$$
$$= \frac{9(P_i^{(n)})^2}{4(P_i^{(n)} + N_i^{(n)})^2} + \frac{x_i^2 N_i^{(n)}}{(P_i^{(n)} + N_i^{(n)})^2}$$
$$\leq \frac{9}{4} + \frac{x_i^2}{A_P^{(n)}}$$
$$\leq \frac{9}{4} + \frac{x_i^2}{P}, \tag{6.64}$$

where we have used $P_i^{(n)} + N_i^{(n)} \geq A_P^{(n)}$ and $A_P^{(n)} \geq P$ (cf. (6.45) and (6.47)). Since $I(\overline{Y}_i^{(n)}|x_i)$ ($i = 1, 2, \cdots, n$) are mutually independent under the conditional distribution $\overline{W}^n(\cdot|\mathbf{x})$, it follows from (6.64) that

$$V_{\mathbf{x}}\left[\frac{1}{n}\sum_{i=1}^n \log \frac{\overline{W}_i(\overline{Y}_i^{(n)}|x_i)}{P_{\tilde{Y}_i}(\overline{Y}_i^{(n)})}\right] \leq \frac{9}{4n} + \frac{\sum_{i=1}^n x_i^2}{n^2 P}$$
$$\leq \frac{13}{4n} \equiv \frac{\sigma_0^2}{n}, \tag{6.65}$$

where we have taken account of input cost constraint (6.48). Chebyshev inequality together with (6.63) and (6.65) leads to

$$\Pr\left\{\frac{1}{n}\sum_{i=1}^n I(\overline{Y}_i^{(n)}|x_i) \geq C_n + \delta \,\bigg|\, \overline{X}^n = \mathbf{x}\right\} \leq \frac{\sigma_0^2}{n\delta^2}, \tag{6.66}$$

where $\delta > 0$ is an arbitrarily small constant and, for simplicity, we have put

$$C_n \equiv \frac{1}{2n}\sum_{i=1}^n \log\left(1 + \frac{P_i^{(n)}}{N_i^{(n)}}\right). \tag{6.67}$$

We notice here that inequality (6.66) holds for all realizations $\mathbf{x}$ of $\overline{X}^n$ with $\overline{\mathbf{X}} = \{\overline{X}^n\}_{n=1}^\infty \in \mathcal{S}_P$. Therefore,

$$\Pr\left\{\frac{1}{n}\sum_{i=1}^n I(\overline{Y}_i^{(n)}|\overline{X}_i^{(n)}) \geq C_n + \delta\right\} \leq \frac{\sigma_0^2}{n\delta^2}. \tag{6.68}$$

Then, from (6.60), (6.62) and (6.68), it follows that

$$\Pr\left\{\frac{1}{n}\log\frac{\overline{W}^n(\overline{Y}^n|\overline{X}^n)}{P_{\tilde{Y}^n}(\overline{Y}^n)} \geq C_n + \delta\right\} \leq \frac{\sigma_0^2}{n\delta^2}.$$

Hence,

$$\lim_{n\to\infty}\Pr\left\{\frac{1}{n}\log\frac{\overline{W}^n(\overline{Y}^n|\overline{X}^n)}{P_{\tilde{Y}^n}(\overline{Y}^n)} \geq C_n + \delta\right\} = 0.$$

Since $\delta > 0$ is arbitrary, we have

$$\begin{aligned}\text{p-}\limsup_{n\to\infty}\frac{1}{n}\log\frac{\overline{W}^n(\overline{Y}^n|\overline{X}^n)}{P_{\tilde{Y}^n}(\overline{Y}^n)} &\leq \limsup_{n\to\infty} C_n \\ &= \lim_{n\to\infty} C_n.\end{aligned}$$

Thus, by (6.55) and (6.59), we have

$$\text{p-}\limsup_{n\to\infty}\frac{1}{n}\log\frac{\overline{W}^n(\overline{Y}^n|\overline{X}^n)}{P_{\overline{Y}^n}(\overline{Y}^n)} \leq \lim_{n\to\infty} C_n,$$

that is,

$$\overline{I}(\overline{\mathbf{X}};\overline{\mathbf{Y}}) \leq \lim_{n\to\infty} C_n.$$

Since $\overline{\mathbf{X}}$ was arbitrary as far as $\overline{\mathbf{X}} \in \mathcal{S}_P$, we conclude that

$$\begin{aligned}\sup_{\overline{\mathbf{X}}\in\mathcal{S}_P}\overline{I}(\overline{\mathbf{X}};\overline{\mathbf{Y}}) &\leq \lim_{n\to\infty} C_n \\ &= \lim_{n\to\infty}\frac{1}{2n}\sum_{i=1}^{n}\log\left(1 + \frac{P_i^{(n)}}{N_i^{(n)}}\right),\end{aligned} \qquad (6.69)$$

which was what to be proven. □

Finally, we have the following lemma.

**Lemma 6.6** The ANWGN channel $\mathbf{W}$ with cost constraint (6.48) satisfies both of assumptions (6.2) and (6.5).

*Proof:*

A simple calculation using (6.36) shows that

$$D(W^n(\,\cdot\,|\mathbf{v})\|W^n(\,\cdot\,|\mathbf{x})) = \frac{1}{2}(\mathbf{v}-\mathbf{x})V_n^{-1}(\mathbf{v}-\mathbf{x})^{\mathrm{T}},$$

where $V_n$ is the covariance matrix as defined in (6.35). Then, the Fisher information matrix $F_n(\mathbf{x})$ is calculated as

$$F_n(\mathbf{x}) = V_n^{-1}.$$

Hence, the norm $||F_n(\mathbf{x})||$ can be expressed as

$$||F_n(\mathbf{x})|| = \frac{1}{\min\left(N_1^{(n)}, N_2^{(n)}, \cdots, N_n^{(n)}\right)}, \qquad (6.70)$$

where $N_1^{(n)}, N_2^{(n)}, \cdots, N_n^{(n)}$ are the eigen values of $V_n$ (cf. (6.37)). On the other hand, from (6.49) in Lemma 6.3 we see that for any positive constant $\delta > 0$ there exists an $n(\delta)$ such that

$$\frac{1}{2n}\sum_{i=1}^{n}\log\left(1 + \frac{P_i^{(n)}}{N_i^{(n)}}\right) \leq C_s(P|\mathbf{W}) + \delta \quad (\forall n \geq n(\delta)). \qquad (6.71)$$

Therefore,

$$\frac{1}{2n}\log\left(1 + \frac{P_i^{(n)}}{N_i^{(n)}}\right) \leq C_s(P|\mathbf{W}) + \delta \quad (\forall n \geq n(\delta);\ \forall i = 1, 2, \cdots, n),$$

from which it follows that

$$\log\left(\frac{A_i^{(n)}}{N_i^{(n)}}\right) \leq 2n(C_s(P|\mathbf{W}) + \delta) \quad (\forall n \geq n(\delta);\ \forall i = 1, 2, \cdots, n). \qquad (6.72)$$

Now, in view of (6.47), a consequence of (6.72) is

$$\log\left(\frac{P}{N_i^{(n)}}\right) \leq 2n(C_s(P|\mathbf{W}) + \delta) \quad (\forall n \geq n(\delta);\ \forall i = 1, 2, \cdots, n). \qquad (6.73)$$

Thus, we have

$$N_i^{(n)} \geq Pe^{-2n(C_s(P|\mathbf{W})+\delta)} \quad (\forall n \geq n(\delta);\ \forall i = 1, 2, \cdots, n). \qquad (6.74)$$

Therefore,

$$\min\left(N_1^{(n)}, N_2^{(n)}, \cdots, N_n^{(n)}\right) \geq Pe^{-2n(C_s(P|\mathbf{W})+\delta)}. \qquad (6.75)$$

Substitution of (6.75) into (6.70) yields

$$||F_n(\mathbf{x})|| \leq \frac{1}{P} e^{2n(C_s(P|\mathbf{W})+\delta)}, \tag{6.76}$$

which obviously satisfies assumption (6.5). Furthermore, it has been shown already in the previous subsection that input cost constraint (6.48) satisfies assumption (6.2). □

Thus, summarizing up Lemmas 6.4 ∼ 6.6 and applying Theorem 6.2 to the ANWGN channel $\mathbf{W}$, we have the following theorem.

**Theorem 6.4** Let $\mathbf{W}$ be the ANWGN channel with input cost constraint (6.48). Then, for all $\varepsilon \geq 0$, $\delta \geq 0$, $\mu \geq 0$, $\lambda \geq 0$ such that

$$\varepsilon < 1, \ \mu + \lambda < 1, \ \delta < 1,$$

it holds that

$$\begin{aligned} C_s(\varepsilon, P|\mathbf{W}) &= D(\mu, \lambda, P|\mathbf{W}) \\ &= S(\delta, P|\mathbf{W}) \\ &= C_s(P|\mathbf{W}) \\ &= \frac{1}{4\pi} \int_{-\pi}^{\pi} \log\left(1 + \frac{f(\lambda)}{g(\lambda)}\right) d\lambda. \end{aligned} \tag{6.77}$$

**Corollary 6.5** Let $\mathbf{W}$ be the ANWGN channel with input cost constraint (6.48). Then, for all $\mu \geq 0$, $\lambda \geq 0$ such that $\mu + \lambda < 1$, it holds that

$$D(\mu, \lambda, P|\mathbf{W}) = C_s(P|\mathbf{W}) = \frac{1}{4\pi} \int_{-\pi}^{\pi} \log\left(1 + \frac{f(\lambda)}{g(\lambda)}\right) d\lambda.$$

**Corollary 6.6** Let $\mathbf{W}$ be the ANWGN channel with input cost constraint (6.48). Then, for all $\delta$ such that $0 \leq \delta < 1$, it holds that

$$S(\delta, P|\mathbf{W}) = C_s(P|\mathbf{W}) = \frac{1}{4\pi} \int_{-\pi}^{\pi} \log\left(1 + \frac{f(\lambda)}{g(\lambda)}\right) d\lambda.$$

# References

[1] R. Ahlswede and G. Dueck, "Identification via channels," *IEEE Transactions on Information Theory*, vol.IT-35, no.1, pp.15-29, 1989

[2] C. E. Shannon, "A mathematical theory of communication," *Bell System Technical Journal*, vol.27, pp.379-423, pp. 623-656, 1948

[3] T.S. Han and S.Verdú, "New results in the theory of identification via channels," *IEEE Transactions on Information Theory*, vol.IT-38, no.1, pp.14-25, 1992

[4] T.S. Han and S. Verdú, "Approximation theory of output statistics," *IEEE Transactions on Information Theory*, vol.IT-39, no.3, pp. 752-772, 1993

[5] S. Verdú and T.S. Han, "A general formula for channel capacity," *IEEE Transactions on Information Theory*, vol.IT-40, no.4, pp.1147-1157, 1994

[6] T. S. Han, *Information-Spectrum Methods in Information Theory*, Baifukan-Press, Tokyo, 1998 (in Japanese)

[7] A. Feinstein, "On the coding theorem and its converse for finite-memory channels," *Information and Control*, vol.2, no.1, pp.25-44, 1959

[8] J. Wolfowitz, "Strong converse of the coding theorem for the general discrete finite-memory channel," *Information and Control*, vol.3, no.1, pp.89-93, 1960

[9] J. Wolfowitz, *Coding Theorems of Information Theory*, 3rd ed., Springer-Verlag, New York, 1978

[10] C.E. Shannon, "Probability of error for optimal codes in a Gaussian channel," *Bell System Tech. J.*, vol.38, no.3, pp. 611-656, 1959

[11] I. Csiszár and J. Körner, *Information Theory: Coding Theorems for Discrete Memoryless Systems*, Academic Press, New York, 1981

[12] T. M. Cover and J. A. Thomas, *Elements of Information Theory*, Wiley, New York, 1991

[13] R. G. Gallager, *Information Theory and Reliable Communication*, John Wiley & Sons, New York, 1968

[14] S. Ihara, *Information Theory for Continuous Systems*, World Scientific, New Jersey, 1993